\documentclass[12pt]{amsart}

\usepackage[margin=1in]{geometry}
\usepackage{amsmath,amssymb,amsthm,mathtools}
\usepackage{amscd}
\usepackage{float}
\usepackage{enumitem}
\usepackage{microtype}
\usepackage{xcolor}

\usepackage{hyperref}
\usepackage{algpseudocode} 
\usepackage{titlesec} 
\titleformat{\section}{\normalfont\Large\bfseries}{\thesection}{0.8em}{}
\titleformat{\subsection}{\normalfont\large\bfseries}{\thesubsection}{0.7em}{}
\titlespacing*{\section}{0pt}{2.2ex plus .3ex minus .2ex}{1.0ex}
\titlespacing*{\subsection}{0pt}{1.7ex plus .2ex minus .2ex}{0.7ex}
\makeatletter
\let\wb@afterheading\@afterheading
\def\@afterheading{\wb@afterheading\@afterindenttrue}
\makeatother

\hypersetup{
  colorlinks=true,
  linkcolor=blue!55!black,
  urlcolor=blue!55!black,
  pdftitle={Computing pro-$p$ Closures of Finitely Generated Subgroups
of Closed Orientable Surface Groups}
}

\theoremstyle{plain}
\newtheorem{theorem}{Theorem}[section]
\newtheorem{proposition}[theorem]{Proposition}
\newtheorem{lemma}[theorem]{Lemma}
\newtheorem{corollary}[theorem]{Corollary}

\theoremstyle{definition}
\newtheorem{definition}[theorem]{Definition}

\newfloat{algorithm}{thp}{lop}[section]
\floatname{algorithm}{Algorithm}

\newtheorem*{mainthm}{Main theorem}

\newcommand{\Cl}{\operatorname{Cl}}
\newcommand{\Core}{\operatorname{Core}}

\newcommand{\id}{\operatorname{id}}

\newcommand{\Fp}{\mathbb F_p}
\newcommand{\Ftwo}{\mathbb F_2}

\newcommand{\cyc}{\operatorname{cyc}}

\begin{document}

\title[Computing pro-$p$ Closures of Surface Groups]{Computing pro-$p$ Closures of Finitely Generated Subgroups
of Closed Orientable Surface Groups}
\author{Zheng Dang}
\address{Department of Mathematics, Soochow University, Suzhou 215006, CHINA}
\email{20244207037@stu.suda.edu.cn}

\author{Jianchun Wu}
\address{Department of Mathematics, Soochow University, Suzhou 215006, CHINA}
\email{wujianchun@suda.edu.cn}

\thanks{$*$The second author is the corresponding author. The authors acknowledge partial support from the National Natural Science Foundation of China (Grant No. 12271385).}

\subjclass[2020]{Primary 20F10; Secondary 20E18, 20F65, 57K20}

\keywords{pro-$p$ topology, pro-$p$ closure, surface group,
core surface, computability}


\begin{abstract}
Let $\Gamma_g$ be the fundamental group of a closed orientable surface of genus
$g\geq2$, and let $p$ be a prime. For every finitely generated subgroup
$H\leq\Gamma_g$ we give an algorithm which, given any finite generating set of $H$ by words in
a standard generating set, outputs a finite set of words in the same alphabet
that generates the pro-$p$ closure of $H$. This extends the
computation of pro-$p$ closures in free groups by Ribes--Zalesskii and
Margolis--Sapir--Weil.
\end{abstract}

\maketitle

\section{Introduction}\label{sec:intro}

A {\it pseudovariety} of finite groups is a class of finite
groups closed under taking subgroups, homomorphic images and finite direct
products. Each such class $\mathbf V$ equips every
group $G$ with a topology, the {\it pro-$\mathbf V$ topology}. The construction
goes back to Hall \cite{hall50}.

Let $p$ be a prime and $\mathbf G_p$ the
pseudovariety of finite $p$-groups; the pro-$\mathbf G_p$ topology is usually called the
{\it pro-$p$ topology}. It has
become a standard finiteness device in combinatorial and geometric group theory,
since finite $p$-group quotients detect the elements of many groups of interest.
For a subgroup
$H\leq G$, one then asks how to describe its closure $\Cl_p(H)$, and whether a
finitely generated $H$ admits an algorithm to compute it.

For a free group the answer is known. A subgroup is encoded by its Stallings
core graph, and the closure is computed from it. This computation is due to
Ribes and Zalesskii \cite{rzk-add}. Margolis,
Sapir and Weil \cite{msw01} gave the algorithmic treatment and settled the
nilpotent case (see also \cite{weil98}). The same questions have been pursued for other pseudovarieties
$\mathbf V$ since then. When $\mathbf V$ is the pseudovariety of finite abelian
groups, of finite abelian groups with bounded exponent, or of
finite metabelian groups, the answers are all positive \cite{msv-abelian}.
Closedness, denseness and membership in the closure are decidable, and a basis
for the closure can be computed whenever it is finitely generated. When
$\mathbf V$ is the pseudovariety of finite supersolvable
groups, denseness is decidable \cite{msv-supersolvable}, and the closure of a
finitely generated subgroup is finitely generated \cite{chen-wu-supersolvable}.

The present paper takes up the next natural case, that of closed orientable
surface groups $\Gamma_g$ of genus $g\geq2$. These groups are residually
$\mathbf G_p$ for every prime $p$ \cite{paris07}. This makes singletons closed,
but does not describe the closure of a general finitely generated subgroup. The
free-group method does not carry over either. The closure of a subgroup of a
free group is read off from its core graph, a one-dimensional object. The
corresponding object for $\Gamma_g$ is two-dimensional in general, so the free-group recipe has no direct analogue.

Our main result fills this gap.

\begin{mainthm}
For every prime $p$ and every finitely generated subgroup $H\leq\Gamma_g$ there
is an algorithm which, given a finite generating set of $H$ by words in a
standard generating set of $\Gamma_g$, outputs a finite generating set of the
pro-$p$ closure $\Cl_p(H)$ by words in the same alphabet.
\end{mainthm}

The computation is geometric. It follows the surface analogue of Stallings's
theory developed by Magee and Puder \cite{mp-core}. Every subgroup of $\Gamma_g$
carries a {\it core surface}, Magee and Puder's analogue of the core graph. Like the core graph, the core surface is an at most two-dimensional complex
obtained from any finite generating set by a folding process, and it is
finite precisely when the subgroup is finitely generated. Write
$\Core(J)$ for the core surface of a subgroup $J$. An inclusion of subgroups
$1\neq H\leq J$ induces a morphism $\Core(H)\to\Core(J)$.

The algorithm works with a descending chain of clopen subgroups
$K_0=\Gamma_g>K_1>K_2>\cdots$ of strictly increasing index, each carrying the
pro-$p$ topology induced from $\Gamma_g$. Density of $H$ in $K_i$ is detected
by the mod-$p$ abelianisation $\sigma_i:K_i\to(\Fp)^{2g(K_i)}$, with $g(K_i)$
the genus of $\Core(K_i)$. If $\sigma_i(H)=(\Fp)^{2g(K_i)}$, then $H$ is dense in
$K_i$ and $\Cl_p(H)=K_i$. Otherwise $K_{i+1}=\sigma_i^{-1}(\sigma_i(H))$ is a
proper clopen surface subgroup of $K_i$ containing $H$.

Each step multiplies $[\Gamma_g:K_i]$ by at least $p$. The core surface of $H$
alone yields a computable bound $\widehat B$, and finitely many steps suffice
to pass it. Once $[\Gamma_g:K_i]>\widehat B$, the growing process of Magee and
Puder \cite{mp-cover}, applied to the morphism $f_i:\Core(H)\to\Core(K_i)$,
factors it through a compact, connected, boundary reduced subsurface with nonempty boundary
$W_{f_i}\subsetneq\Core(K_i)$:
$$
 \Core(H)\xrightarrow{\ f_i\ }W_{f_i}\xrightarrow{\ \iota_{f_i}\ }\Core(K_i).
$$
Cutting $\Core(K_i)$ along an essential boundary curve of $W_{f_i}$ then yields
a free subgroup $L$ of $K_i$. We show that $H\leq\pi_1(W_{f_i})\leq L$ and that
$L$ is pro-$p$-closed in $\Gamma_g$. The pro-$p$ topology induced on $L$ from
$\Gamma_g$ then agrees with its intrinsic pro-$p$ topology. The closure of $H$
in $L$ is therefore $\Cl_p(H)$. To compute it one writes the given words in a
free basis of $L$ and applies the free-group algorithm of \cite{rzk-add,msw01}.
The main theorem follows.

The paper is organised as follows: Section~\ref{sec:prelim} sets up the
pro-$\mathbf V$ topology and core surfaces, Section~\ref{sec:main} the density
criterion, the pro-$p$-closedness of subgroups generated by standard subsets,
and the reduction to free subgroups, Section~\ref{sec:toolkit} the algorithms
on surfaces, and Section~\ref{sec:algorithm} the algorithm and its correctness.

\section{Preliminaries}\label{sec:prelim}

\subsection{The pro-$\mathbf V$ topology}\label{sec:basicdef}

\begin{definition}\label{def:topology}
Let $\mathbf V$ be a pseudovariety of finite groups and $G$ a group. The collection
$$\{gN\mid g\in G,\ N\trianglelefteq G,\ G/N\in\mathbf V\}$$
is a basis of a topology on $G$, called the {\it pro-$\mathbf V$ topology};
topological properties are prefixed with {\it pro-$\mathbf V$}. For
$\mathbf V=\mathbf G_p$ this is the pro-$p$ topology of the introduction.

For $K\leq G$, the pro-$\mathbf V$
topology on $K$ may differ from the subspace topology induced from $G$. If they coincide, then we say the pro-$\mathbf V$
topology of $K$ is {\it induced from} $G$.

We write $\Cl_p(H)$ for the closure of $H\subseteq G$, and, for $H\leq K\leq G$,
$\Cl_p(H,K)$ for the closure of $H$ inside $K$ with respect to the pro-$p$
topology on $K$. We say that $H$ is
{\it pro-$p$ dense} in $K$ if $\Cl_p(H,K)=K$.
\end{definition}

\begin{definition}\label{def:residual}
A homomorphism $\phi:G\to V$ with $V\in\mathbf V$ {\it separates} $x,y\in G$
if $\phi(x)\neq\phi(y)$. If every pair of distinct elements is separated by some
member of $\mathbf V$, then $G$ is {\it residually $\mathbf V$}.
\end{definition}

\begin{definition}\label{def:retract}
A subgroup $H\leq G$ is a {\it retract} if there is a homomorphism
$\phi:G\to H$ with $\phi\circ\jmath=\id_H$, where $\jmath:H\to G$ is the
inclusion.
\end{definition}

\begin{proposition}[{\cite[Proposition~1.5]{msw01}}]\label{prop:retract-topology}
If $H$ is a retract of $G$, then the pro-$\mathbf V$ topology of $H$ is induced
from $G$.
\end{proposition}

\begin{proposition}[{\cite[Corollary~1.3]{msw01}}]\label{prop:retract-closed}
Let $G$ be residually $\mathbf V$ with its pro-$\mathbf V$ topology, and let
$H$ be a retract of $G$. Then $H$ is closed in $G$.
\end{proposition}

\begin{definition}\label{def:extclosed}
A pseudovariety $\mathbf V$ is {\it extension closed} if $N\trianglelefteq G$,
$N,G/N\in\mathbf V$ imply $G\in\mathbf V$. 
\end{definition}

\begin{proposition}[{\cite[Proposition~1.6]{msw01}}]\label{prop:clopen-topology}
Suppose $\mathbf V$ is extension closed, and let $H$ be a pro-$\mathbf V$ clopen subgroup of $G$. Then the pro-$\mathbf V$ topology of
$H$ is induced from $G$.
\end{proposition}

For $p$ a
prime we write $\mathbf{Ab}_p$ for the pseudovariety of finite abelian $p$-groups.

\begin{lemma}\label{lem:density-iff-abelian}
Let $G$ be any group, and let $H\leq G$. Then:
\[
 H\text{ is pro-$p$ dense}\quad\Longleftrightarrow\quad
 H\text{ is pro-$\mathbf{Ab}_p$ dense}.
\]
\end{lemma}
\begin{proof}
$(\Rightarrow)$. Since $\mathbf{Ab}_p\subseteq\mathbf G_p$, every
pro-$\mathbf{Ab}_p$ open set is pro-$p$ open; equivalently, the
pro-$\mathbf{Ab}_p$ topology is coarser than the pro-$p$ topology. A subset dense
in a finer topology is dense in a coarser one, so pro-$p$ density implies
pro-$\mathbf{Ab}_p$ density.

$(\Leftarrow)$. Suppose $H$ is not pro-$p$ dense. Then $H$ lies in a proper subgroup that is open in the pro-$p$ topology. This subgroup contains some $N\trianglelefteq
G$ with $G/N$ a finite $p$-group, and the quotient map $\phi:G\to P:=G/N$ then
sends $H$ into a proper subgroup of $P$. A proper subgroup of a finite $p$-group lies in an
index-$p$ normal subgroup, so composing $\phi$ with the quotient
$P\twoheadrightarrow\Fp$ gives an epimorphism $\psi:G\to\Fp$ with
$\psi(H)=0$. Since $\Fp\in\mathbf{Ab}_p$, the kernel $\ker\psi$ is open,
hence clopen, in the pro-$\mathbf{Ab}_p$ topology, and it is proper and contains
$H$. The closure of $H$ is contained in $\ker\psi$, which is proper, so $H$ is
not pro-$\mathbf{Ab}_p$ dense.
\end{proof}

\subsection{Core surfaces}\label{sec:cores}
For a closed connected orientable surface $\Sigma_g$ of genus $g\geq2$ with a
basepoint $o$, let
\[
 \Gamma_g=\pi_1(\Sigma_g,o)=\Bigl\langle a_1,b_1,\dots,a_g,b_g
 \ \Big|\ [a_1,b_1]\cdots[a_g,b_g]=1\Bigr\rangle
\]
be its fundamental group. Throughout,
$[a,b]=aba^{-1}b^{-1}$, and we call
$\{a_1,b_1,\dots,a_g,b_g\}$ a {\it standard generating set}. $\Sigma_g$ carries a CW-structure with one vertex $o$, the $2g$ directed edges $a_1,b_1,\dots,a_g,b_g$, and one $4g$-gon attached along the word
$[a_1,b_1]\cdots[a_g,b_g]$. This CW-structure pulls back to every covering
surface of $\Sigma_g$. We
measure the length of a path in the $1$-skeleton by its number of edges.

A subgroup of a free group is encoded by its \emph{Stallings core graph}. This
graph is the subgraph of the covering graph obtained by pruning the hanging
trees. It is finite exactly when the subgroup is finitely generated, and it is
built from any finite generating set by Stallings folds
\cite{kapovich-mi,stallings}. Magee and Puder introduced the surface group
analogue of this graph, and the definition below follows them. The covering
surfaces corresponding to conjugate subgroups are isomorphic, via an isomorphism
induced by a deck transformation of the universal cover, so the core surface
defined below depends only on the conjugacy class of the subgroup.

\begin{definition}[{\cite[Definition~1.1]{mp-core}}]\label{def:core}
For a subgroup $J\leq\Gamma_g$ let $\Upsilon_J$ be the
covering surface corresponding to $J$ (with the pulled back CW-structure). The
{\it core surface} $\Core(J)$ of $J$ is the subcomplex of $\Upsilon_J$ obtained as
follows. Take the union of all shortest representatives in $\Upsilon_J^{(1)}$ of the free
homotopy classes of essential closed curves, and adjoin every connected
component of the complement of this union that contains only finitely many $2$-cells. For
completeness, the core surface of the trivial subgroup is the $0$-dimensional
complex consisting of a single vertex mapped to $o$.
\end{definition}

To compute with core surfaces we need a combinatorial model. Both the
construction of $\Core(J)$ from a finite generating set and the algorithms of
Section~\ref{sec:toolkit} are carried out on the tiled surfaces below, using
their combinatorial data. Throughout, a {\it signed edge} is one of the two
{\it orientations} of an edge. We write $e^{\pm}$ for the two signed edges of an
edge $e$, and $\overline{e^{\pm}}=e^{\mp}$ for the reverse signed edge. A {\it
direction} on $e$ is the choice of one signed edge, $e^{+}$. A {\it path} in a
$1$-complex is a sequence of signed edges in which the terminal vertex of each
term is the initial vertex of the next. A path is {\it closed} if it starts and
ends at the same vertex. It {\it reads} the word obtained from the labels of its
edges, with exponent $+1$ for signed edges traversed along the direction and
$-1$ otherwise.

\begin{definition}[{\cite[Proposition~3.3]{mp-core}}]\label{def:tiled}
A ($g$-){\it tiled surface} is an at most two--dimensional CW-complex in which
every edge carries a {\it direction} and a label in $\{a_1,b_1,\dots,a_g,b_g\}$ such that
\begin{enumerate}[label=\textup{(P\arabic*)}]
\item at each vertex there is at most one incoming and one outgoing edge of each
      label;
\item every path in the $1$-skeleton reading a word equal to $1\in\Gamma_g$ is
      closed;
\item every $2$-cell is a $4g$-gon attached along a closed path reading the
      relation $[a_1,b_1]\cdots[a_g,b_g]$, and every such closed path is the
      boundary of at most one such $4g$-gon.
\end{enumerate}
\end{definition}
Equivalently, a tiled surface $Y$ can be realised as a subcomplex of a (not
necessarily connected) covering surface of $\Sigma_g$ (see \cite[Definition~3.1]{mp-core}).
The directions and labels are pulled back from $\Sigma_g$ and serve to
describe the restricted covering map $p:Y\to\Sigma_g$.
A {\it morphism} between tiled surfaces is a map of their CW-complexes sending
$i$-cells to $i$-cells for $i=0,1,2$ and respecting the directions and labels of
edges \cite[p.~46]{mp-core}. If a tiled surface $Y$ is contained in another tiled
surface $W$, we call $Y$ a tiled \emph{subsurface} of $W$.

\begin{definition}\label{def:celldata}
For a compact tiled surface $Y$, its {\it combinatorial data} are the following:
\begin{enumerate}[label=\textup{(\roman*)}]
\item the vertex set $V$, edge set $E$ and $2$-cells (each a $4g$-gon), with
      the incidence of vertices on edges;
\item for each edge $e$ its two signed edges $e^{+},e^{-}$, with $e^{+}$ the one
      in the direction assigned in Definition~\ref{def:tiled} and
      $\overline{e^{\pm}}=e^{\mp}$ the reverse;
\item at each vertex $v$ a cyclic order $\cyc(v)$ of the signed edges {\it
      leaving} $v$, inherited from the vertex of $\Sigma_g$: it is the cyclic
      subsequence, consisting of the signed edges present at $v$, of
      (outgoing $a_1$, incoming $b_1$, incoming $a_1$, outgoing $b_1$,
      outgoing $a_2$, incoming $b_2$, incoming $a_2$, outgoing $b_2$, $\dots$).
      Gluing the $2$-cells along matching labels determines $\cyc(v)$
      symbolically from the labels alone;
\item for each $2$-cell $O$ (a $4g$-gon), $P(O)$ is the closed path that goes
      once around the boundary of $O$ counterclockwise. We will use $P(O)$ as
      the attaching map of $O$.
\end{enumerate}
\end{definition}


\begin{definition}[{\cite[Definition~3.2]{mp-core}}]\label{def:thick}
Let $Y$ be a tiled surface, realised as a subcomplex of a covering surface
$\Upsilon$ of $\Sigma_g$. A small closed regular
neighbourhood $N(Y)$ of $Y$ in $\Upsilon$ is called the {\it thick version} (or {\it thickening}) of $Y$. It is a surface, possibly with
boundary. It does not depend on the ambient surface $\Upsilon$ and can be
recovered from the combinatorial data of $Y$ (see
\cite[Proposition~3.4]{mp-core}).
A {\it boundary cycle} of $Y$ is a closed path in $Y^{(1)}$ corresponding to a
boundary component of $N(Y)$, with the orientation of
\cite[Section~3]{mp-core}.
\end{definition}

Below we use two boundary conditions on tiled surfaces. A tiled surface $Y$ is \emph{boundary
reduced} if no boundary cycle of $Y$ contains a long block or a long chain, and
\emph{strongly boundary reduced} if no boundary cycle of $Y$ contains a
half-block or a half-chain; see \cite[Sections~3--4]{mp-core} for the definitions.
Strongly boundary reduced implies boundary reduced.

\begin{proposition}[{\cite[Corollary~4.11]{mp-core}}]\label{prop:pi1-injective}
Every morphism from a connected boundary reduced tiled surface is
$\pi_1$-injective.
\end{proposition}

We collect the properties of core surfaces used in this paper.

\begin{proposition}[{\cite[Proposition~5.3, 5.8, Lemma~5.4, Theorem~5.10]{mp-core}}]\label{prop:mp-basic}
Let $\{1\}\neq J\leq\Gamma_g$ and $\Upsilon_J$ be the
covering surface corresponding to $J$.
\begin{enumerate}[label=\textup{(\arabic*)}]
\item $\Core(J)$ is connected and strongly boundary reduced, and the embedding
$\Core(J)\hookrightarrow\Upsilon_J$ induces an isomorphism on fundamental groups.
\item $\Core(J)$ is compact if and only if $J$ is finitely generated.
\item (Naturality) If $1\ne H\leq J$, the covering map restricts to a morphism
of tiled surfaces $\Core(H)\to\Core(J)$.
\item (Construction) Given any finite generating set of $J$, the core
surface $\Core(J)$ is produced from the labelled wedge on the generators by a
finite folding process.
\end{enumerate}
\end{proposition}

\section{The pro-$p$ closure of finitely generated subgroups}\label{sec:main}

For a subgroup $H\leq\Gamma_g$, the mod-$p$ abelianisation gives a simple test for pro-$p$ density, that is, for whether $\Cl_p(H)=\Gamma_g$.

\begin{proposition}\label{prop:density-equiv}
Let
\[
 \sigma:\Gamma_g\longrightarrow \Gamma_g^{\mathrm{ab}}\cong\mathbb Z^{2g}
 \xrightarrow{\ \mathrm{mod}\ p\ }(\Fp)^{2g}
\]
be the composite of abelianisation and reduction modulo $p$. For a  subgroup $H\leq\Gamma_g$, the following are equivalent:
\begin{enumerate}[label=\textup{(\arabic*)}]
\item $H$ is pro-$p$ dense in $\Gamma_g$;
\item $H$ is pro-$\mathbf{Ab}_p$ dense in $\Gamma_g$;
\item $\sigma(H)=(\Fp)^{2g}$;
\item $\sigma^{-1}(\sigma(H))=\Gamma_g$.
\end{enumerate}
\end{proposition}
\begin{proof}
By Lemma~\ref{lem:density-iff-abelian}, (1) and (2) are equivalent. Since $\sigma$ is surjective, (3) and (4) are equivalent.

(2)$\Rightarrow$(3): Suppose $\sigma(H)$ is a
proper subspace of $(\Fp)^{2g}$. Since this space is finite-dimensional, some
nonzero linear functional $\lambda:(\Fp)^{2g}\to\Fp$ vanishes on $\sigma(H)$, and
then $\lambda\circ\sigma:\Gamma_g\to\Fp$ is an epimorphism whose kernel
contains $H$. This kernel is clopen in the pro-$\mathbf{Ab}_p$ topology and proper,
so $H$ is not pro-$\mathbf{Ab}_p$ dense.

(3)$\Rightarrow$(2): Conversely, suppose $H$ lies in a proper open subgroup $K$
of the pro-$\mathbf{Ab}_p$ topology. Then $K$ contains the kernel of an
epimorphism $\psi:\Gamma_g\to A$ onto a finite abelian $p$-group $A$. Since
$K$ is proper, the image $\psi(H)\subseteq \psi(K)$ is a proper subgroup of $A$. A proper subgroup
of a finite $p$-group lies in an index-$p$ subgroup, so some quotient
$A\twoheadrightarrow\Fp$ vanishes on $\psi(H)$. Composing with $\psi$ yields an
epimorphism $\chi:\Gamma_g\to\Fp$ with $\chi(H)=0$. Its target is abelian of
exponent $p$, so $\chi$ factors as $\chi=\lambda\circ\sigma$ for a nonzero
linear functional $\lambda:(\Fp)^{2g}\to\Fp$. Since $\chi(H)=0$, $\lambda$
vanishes on $\sigma(H)$, making $\sigma(H)$ a proper subspace of $(\Fp)^{2g}$.
This contradicts (3).
\end{proof}

For some special subgroups $H\leq\Gamma_g$, the pro-$p$ topology of $H$ is induced from $\Gamma_g$.

\begin{lemma}\label{lem:small-retract}
Let $S$ be a subset of a standard generating set $\{a_1,b_1,\dots,a_g,b_g\}$ of $\Gamma_g$ with $|S|\le g$. Then $H=\langle S\rangle$ is a retract of
$\Gamma_g$.
\end{lemma}
\begin{proof}
The set $S$ is a proper subset of a standard generating set. By the Freiheitssatz \cite{lyndon-schupp}, $S$ freely generates $H$. Partition the indices $\{1,\dots,g\}$ into
\[
 I_2=\{i: a_i,b_i\in S\},\quad
 I_1=\{i:\text{ exactly one of }a_i,b_i\in S\},\quad
 I_0=\{i: a_i,b_i\notin S\}.
\]
Then
$|S|=2|I_2|+|I_1|$ while $g=|I_2|+|I_1|+|I_0|$, so the hypothesis $|S|\le g$ gives
$|I_0|\ge|I_2|$.

We construct an injection $\sigma:I_2\to I_0$ recursively and record the pairs
it outputs in the order produced,
\[
 (i_1,\sigma(i_1)),\ (i_2,\sigma(i_2)),\ \dots,\ (i_{|I_2|},\sigma(i_{|I_2|})).
\]
Run the recursion on temporary copies $J_2:=I_2$ and $J_0:=I_0$, so that $I_2$
and $I_0$ themselves are left unchanged. Applied to the current sets $J_2$ and
$J_0$, a step outputs a single pair as follows:
\begin{enumerate}[label=\textup{(\roman*)}]
\item If some element of $J_0$ exceeds $\max J_2$, output
$\bigl(\max J_2,\ \min\{j\in J_0:j>\max J_2\}\bigr)$.
\item Otherwise every element of $J_0$ is less than $\max J_2$; output
$\bigl(i,\ \max J_0\bigr)$ with $i=\min\{i'\in J_2:i'>\max J_0\}$.
\end{enumerate}
Delete the two elements of the output pair from $J_2$ and $J_0$ respectively, and repeat until
$J_2=\varnothing$. Each step matches a new forced index to a new available index, so $\sigma$
is a well-defined injection. Since $|J_0|\ge|J_2|$ throughout, an available partner exists at every step
and the procedure terminates.

Define $f$ on the standard generators by
\[
 f(x)=
 \begin{cases}
  x, & x\in S,\\
  1, & x\in\{a_i,b_i\},\ i\in I_1,\ x\notin S,\\
  b_i, & x=a_{\sigma(i)},\ i\in I_2,\\
  a_i, & x=b_{\sigma(i)},\ i\in I_2,\\
  1, & x\in\{a_i,b_i\},\ i\in I_0\setminus\sigma(I_2)
 \end{cases}
\]
The group $\Gamma_g$ is presented by
$\langle a_1,b_1,\dots,a_g,b_g\mid[a_1,b_1]\cdots[a_g,b_g]\rangle$. Hence $f$ extends to a
homomorphism $\Gamma_g\to H$ precisely when $f$ sends the relator
$[a_1,b_1]\cdots[a_g,b_g]$ to $1$ in $H$, that is, when
\begin{equation}\label{eq:relimg}
 f\bigl([a_1,b_1]\cdots[a_g,b_g]\bigr)
 =\prod_{i=1}^{g}\bigl[f(a_i),f(b_i)\bigr]
\end{equation}
equals $1$.

 The factor of \eqref{eq:relimg} at
a position $i\in I_2$ is $[a_i,b_i]$, the factor at position $\sigma(i)$ is its inverse, and
every position in $I_1\cup(I_0\setminus\sigma(I_2))$ contributes $1$. Cancel the pairs in the
order $i_1,i_2,\dots,i_{|I_2|}$ in which they were produced. When the pair
$i_k$ is reached, its two endpoints carry $[a_{i_k},b_{i_k}]$ and its inverse, so it remains to
account for what lies strictly between them. At the step that produced $i_k$, in
case~(i) $i_k$ is $\max J_2$ and $\sigma(i_k)$ is the least element of $J_0$
above it; in case~(ii) $\sigma(i_k)=\max J_0$ and $i_k$ is the least element of
$J_2$ above $\max J_0$. In both cases a position strictly between the endpoints
is either trivial or a position of a pair matched at an earlier step (hence
already suppressed). Thus the subword between $\min(i_k,\sigma(i_k))$ and
$\max(i_k,\sigma(i_k))$ consists of two mutually inverse factors,
$[a_{i_k},b_{i_k}]$ and its inverse, in one order or the other, and is deleted.
Proceeding through
$k=1,\dots,|I_2|$ deletes every pair, so \eqref{eq:relimg} is the identity in $H$. Hence $f$
extends to a homomorphism of $\Gamma_g$, and since $f$ fixes $S$ it restricts to the identity on $H$; $f$ is a retraction of $\Gamma_g$ onto $H$.
\end{proof}

\begin{proposition}\label{prop:standard-retract}
Let $S$ be a proper subset of a standard generating set $\{a_1,b_1,\dots,a_g,b_g\}$ of $\Gamma_g$.
Then $H=\langle S\rangle$ is pro-$p$-closed in $\Gamma_g$, and the pro-$p$
topology of $H$ is induced from $\Gamma_g$.
\end{proposition}
\begin{proof}
Let $c$ be a generator not in $S$. Re-indexing the standard generating set
cyclically, that is, replacing
$(a_i,b_i)$ by $(a_{i+1},b_{i+1})$ with indices modulo $g$, turns the relator
into a cyclic permutation of itself, so the re-indexed set is again a standard
generating set. Thus we may assume that $c$ lies in the
last handle; then $c$ is $b_g$ or $a_g$. Since $S$ is a proper subset of the
$2g$ standard generators, $|S|\le 2g-1$. Choose $e\ge1$ large enough that
$g':=1+(g-1)p^e\ge 2g-1$.

If $c=b_g$, let $\varphi:\Gamma_g\to\mathbb Z$ send $b_g$ to $1$ and vanish on
the other standard generators, and let $K'$ be the kernel of the epimorphism
$\Gamma_g\xrightarrow{\ \varphi\ }\mathbb Z\xrightarrow{\ \mathrm{mod}\,p^e\ }\mathbb Z_{p^e}$;
then $K'$ is a clopen subgroup of index $p^e$, hence a genus-$g'$ surface
subgroup. Reidemeister--Schreier with the transversal
$\{1,b_g,\dots,b_g^{p^e-1}\}$, combined with the relation
$[a_1,b_1]\cdots[a_g,b_g]=1$, shows that $K'$ admits the standard generating
set
\[
 S'=\{ a_g,\,b_g^{p^e}\}\cup\{ b_g^k a_i b_g^{-k},\,
 b_g^k b_i b_g^{-k}\mid 1\le i<g,\ 0\le k<p^e\},
\]
which has $2(g-1)p^e+2=2g'$ elements, paired as $(a_g,b_g^{p^e})$ and
$(b_g^k a_i b_g^{-k},\,b_g^k b_i b_g^{-k})_{1\le i<g}$, with relator
\[
 \prod_{i=1}^{g-1}[\,b_g^{p^e-1}a_i b_g^{-(p^e-1)},\,
   b_g^{p^e-1}b_i b_g^{-(p^e-1)}]\;\cdots\;
 \prod_{i=1}^{g-1}[\,b_g a_i b_g^{-1},\,
   b_g b_i b_g^{-1}]
 \prod_{i=1}^{g-1}[\,a_i,\,b_i]\;[a_g,\,b_g^{p^e}]=1.
\]

If $c=a_g$, the same construction with $a_g$ in place of $b_g$ applies. Let
$K'$ be the kernel of the epimorphism $\Gamma_g\to\mathbb Z_{p^e}$ sending
$a_g$ to $1$ and all other standard generators to $0$, again clopen of index
$p^e$ and genus $g'$, and
use the transversal $\{1,a_g,\dots,a_g^{p^e-1}\}$. Reidemeister--Schreier
gives the standard generating set
\[
 S'=\{ a_g^{p^e},\,b_g\}\cup\{ a_g^k a_i a_g^{-k},\,
 a_g^k b_i a_g^{-k}\mid 1\le i<g,\ 0\le k<p^e\},
\]
of $2g'$ elements, paired as $(a_g^{p^e},b_g)$ and
$(a_g^k a_i a_g^{-k},\,a_g^k b_i a_g^{-k})_{1\le i<g}$, with relator
\[
 \prod_{i=1}^{g-1}[\,a_i,\,b_i]\;
 \prod_{i=1}^{g-1}[\,a_g a_i a_g^{-1},\,a_g b_i a_g^{-1}]\;\cdots\;
 \prod_{i=1}^{g-1}[\,a_g^{p^e-1}a_i a_g^{-(p^e-1)},\,
   a_g^{p^e-1}b_i a_g^{-(p^e-1)}]
 \;[a_g^{p^e},b_g]=1.
\]

In both cases $S\subseteq S'$: the transversal is taken along the missing
generator $c$, and $S'$ contains every standard generator except $c$. Since
$|S|\le 2g-1\le g'$, applying Lemma~\ref{lem:small-retract} to the genus-$g'$
surface group $K'$ shows that $H=\langle S\rangle$ is a retract of $K'$. 

Since $\mathbf G_p$
is extension closed and $K'$ is clopen in $\Gamma_g$, the pro-$p$ topology of $K'$ is
induced from $\Gamma_g$ by Proposition~\ref{prop:clopen-topology}. Because $H$ is a retract of $K'$, the pro-$p$ topology of $H$ is
induced from $K'$ by Proposition~\ref{prop:retract-topology}, and hence from $\Gamma_g$. The group $K'$ is residually $\mathbf G_p$
(see \cite{paris07}), so $H$ is
closed in $K'$ by Proposition~\ref{prop:retract-closed}. Finally, $K'$ is closed in
$\Gamma_g$, and therefore $H$ is closed in $\Gamma_g$ too.
\end{proof}

\begin{lemma}\label{lem:automorphism-invariance}
Let $\varphi\colon G\to G$ be an automorphism and $A\le G$. Then $A$ is
pro-$p$-closed in $G$, and its pro-$p$ topology is induced from $G$, if and only
if $\varphi(A)\le G$ has the same two properties.
\end{lemma}
\begin{proof}
The normal subgroups of $G$ of $p$-power index form a neighbourhood base at
the identity for the pro-$p$ topology, and $\varphi$ permutes this collection,
since $G/\varphi(N)\cong G/N$. Hence $\varphi$ is a homeomorphism of the
pro-$p$ topology of $G$, so $A$ is pro-$p$-closed in $G$ if and only if
$\varphi(A)$ is. For the intrinsic topologies, $\varphi$ restricts to a
bijection between the subgroups of $A$ and those of $\varphi(A)$. It carries
$M\trianglelefteq A$ with $A/M$ a finite $p$-group to $\varphi(M)
\trianglelefteq\varphi(A)$ with $\varphi(A)/\varphi(M)\cong A/M$ a finite
$p$-group. These two collections form neighbourhood bases at the identity for
the intrinsic pro-$p$ topologies of $A$ and $\varphi(A)$, so $\varphi$
identifies the intrinsic pro-$p$ topology of $A$ with that of $\varphi(A)$.
Since $\varphi$ is a homeomorphism of $G$, it also identifies the topologies
induced on $A$ and $\varphi(A)$ from $G$. Hence the two properties hold for
$A$ if and only if they hold for $\varphi(A)$.
\end{proof}

\begin{corollary}\label{cor:automorphism-closed}
Let $\varphi\in\operatorname{Aut}(\Gamma_g)$ and $L\leq\Gamma_g$. If $\varphi(L)=\langle S\rangle$ for a proper subset $S$ of a standard generating set of $\Gamma_g$, then $L$ is pro-$p$-closed in $\Gamma_g$ with its pro-$p$ topology induced from $\Gamma_g$.
\end{corollary}
\begin{proof}
By Proposition~\ref{prop:standard-retract}, $\langle S\rangle$ is pro-$p$-closed in $\Gamma_g$ with its pro-$p$ topology induced from $\Gamma_g$. Since $L=\varphi^{-1}(\langle S\rangle)$, Lemma~\ref{lem:automorphism-invariance} transfers these two properties to $L$.
\end{proof}

The following reduces the computation of the closure to the free-group case.

\begin{proposition}\label{prop:free-reduction}
Let $H\leq L\leq K\leq\Gamma_g$, where $K$ is a pro-$p$ clopen subgroup and
$L$ is a free subgroup of $K$ such that $L$ is pro-$p$-closed in $K$ and its
pro-$p$ topology is induced from $K$. Then
\[
 \Cl_p(H,L)=\Cl_p(H,K)=\Cl_p(H,\Gamma_g).
\]
\end{proposition}
\begin{proof}
We use the following observation for $H\le G_1\le G_2$. If the pro-$p$
topology of $G_1$ is induced from $G_2$, then
\[ \Cl_p(H,G_1)=\Cl_p(H,G_2)\cap G_1. \]
If moreover $G_1$ is pro-$p$-closed in $G_2$, then $\Cl_p(H,G_2)\subseteq G_1$,
so the two closures coincide.

Apply this to $H\le L\le K$. By hypothesis the pro-$p$ topology of $L$ is
induced from $K$, and $L$ is pro-$p$-closed in $K$, so $\Cl_p(H,L)=\Cl_p(H,K)$.
Apply it to $H\le K\le\Gamma_g$. Since $K$ is clopen in $\Gamma_g$,
Proposition~\ref{prop:clopen-topology} shows the pro-$p$ topology of $K$ is
induced from $\Gamma_g$; being clopen, $K$ is pro-$p$-closed in $\Gamma_g$, so
the two closures coincide. The conclusion follows.
\end{proof}

\section{Algorithms on surfaces}\label{sec:toolkit}

This section develops the surface algorithms used in
Section~\ref{sec:algorithm}. Section~\ref{sec:freegen} computes a free basis
of the fundamental group of a compact surface with boundary.
Section~\ref{sec:boundary-extract} produces, for a boundary reduced subsurface
of a closed tiled surface, a free subgroup that is pro-$p$-closed with its
topology induced and has a computable free basis.

\subsection{A free basis of a surface with boundary}\label{sec:freegen}

Let $Y$ be a compact connected orientable surface with nonempty
boundary, given as a finite cell complex whose $2$-cell boundary paths meet each
signed edge at most once. Orient $Y$. Each $2$-cell inherits an orientation
from $Y$, and $P(O)$ is the closed path going once around the boundary of $O$
in that orientation, so two $2$-cells sharing an edge traverse it in opposite
directions. Choose a direction on each $1$-cell, fix a basepoint, and choose a
spanning tree $T$ of the $1$-skeleton. Collapsing $T$ turns $Y^{(1)}$ into a wedge of circles, one for
each $1$-cell $e\notin T$. The loop corresponding to $e$ is a generator $x_e$,
and we set $x_e=1$ for $e\in T$. Reading the boundary path $P(O)$ of a $2$-cell
$O$ once, with $x_e$ for a crossing along $e$ and $x_e^{-1}$ for one against
it, gives a relator $r_O$, so
\[
 \pi_1(Y)\cong\langle\,x_e\ (e\notin T)\mid r_O\ (O\in Y^{(2)})\,\rangle .
\]

The dual graph $D(Y)$ has the $2$-cells $O^*$ as vertices, together with
one vertex $\partial^*$ for the boundary. Each $1$-cell $e$ contributes a dual
edge $e^*$, defined as follows according to where its two signed edges occur.
Suppose a signed edge of $e$ occurs in $P(O)$ and the other signed edge of $e$
occurs in $P(O')$. Then $e^*$ joins $O^*$ to $O'^*$. Suppose a signed edge of $e$ occurs in
$P(O)$ and the other occurs in no $P(O')$. Then $e^*$ joins $O^*$ to
$\partial^*$. Two remaining cases give a loop. If both signed edges of $e$ occur
only in $P(O)$, then $e^*$ is a loop at $O^*$. If neither signed edge of $e$
occurs in any $P(O')$, then $e^*$ is a loop at $\partial^*$. Deleting the dual edges $e^*$ with $e\in T$
leaves $D_T(Y)$. Collapsing $T$ to a point gives a complex $X$ homotopy
equivalent to $Y$, and its $1$-cells are exactly the non-tree ones. The collapse leaves each $1$-cell $e\notin T$ and each $2$-cell $O$ intact,
and the signed edges of $e$ occur in $P(O)$ exactly as before. Hence $D_T(Y)$ is the
dual graph of $X$. We show that
$D_T(Y)$ is connected. Since $Y$ is a connected surface with nonempty boundary,
$H_2(Y;\Ftwo)=0$, and hence $H_2(X;\Ftwo)=0$. Suppose that $D_T(Y)$ has a
component $\mathcal C$ not containing the boundary vertex $\partial^*$, and let
\[
 z=\sum_{O^*\in\mathcal C}O\ \in C_2(X;\Ftwo).
\]
Then $z\neq0$. A $1$-cell of $X$ is an edge $e$ outside $T$. The
coefficient of $e$ in $\partial_2z$ is the number, modulo $2$, of occurrences of
the signed edges of $e$ in the paths $P(O)$ with $O^*\in\mathcal C$. Suppose
$e^*$ joins two $2$-cells. Then these two cells both belong to $\mathcal C$ or
both do not, and the two occurrences cancel. Suppose $e^*$ is a loop at a
$2$-cell $O^*$. Then both occurrences lie in $P(O)$ and cancel. Suppose $e^*$
joins a $2$-cell to $\partial^*$, or is a loop at $\partial^*$. Then no $2$-cell
of $\mathcal C$ contributes an occurrence. In every case the coefficient of $e$
in $\partial_2z$ is $0$. Thus $\partial_2z=0$. Since $X$
has no $3$-cells, $H_2(X;\Ftwo)=\ker\partial_2$, so $z$ represents a nonzero
class in $H_2(X;\Ftwo)$, a contradiction. Hence $D_T(Y)$ is connected.

Choose a spanning tree $\mathbb T$ of $D_T(Y)$,
rooted at $\partial^*$ and oriented away from it. Every non-root vertex $O^*$
has a unique parent edge, a dual edge $e^*$ for some $1$-cell $e$ of $Y$. The
edges of $T$ were deleted, so $e\notin T$. Write $e$ as $e_O$. The parent
edge joins $O^*$ to $\partial^*$ or to the vertex $O_0^*$ of another $2$-cell
$O_0$, which we call the {\it parent cell} of $O$. Distinct non-root vertices
have distinct parent edges, so $O\mapsto e_O$ injects the non-root $2$-cells
into the $1$-cells outside $T$, with image the edge set of $\mathbb T$.

Since $e_O^*$ is an edge of the tree $\mathbb T$, it is not a loop. Hence
exactly one of the two signed edges of $e_O$ occurs in $P(O)$, and since $P(O)$
traverses $e_O$ exactly once, the generator
$x_{e_O}$ occurs in $r_O$ exactly once, with exponent $\pm1$. Moreover, for
$O'\neq O$ the generator $x_{e_{O'}}$ occurs in $r_O$ only when $O$ is the
parent cell of $O'$ in $\mathbb T$, since the two sides of $e_{O'}$ are $O'$
and the cell at the other end of the parent edge, and when that end is
$\partial^*$ the other signed edge of $e_{O'}$ occurs in no relator.

We eliminate the generators $x_{e_O}$ in order of decreasing distance
from $\partial^*$ in $\mathbb T$. Thus each $2$-cell is treated only after all
its children. Let $O$ be a $2$-cell, and suppose all its children have already
been treated. Let $r_O$ be its relator at that stage. By the preceding paragraph
$r_O$ is a word in $x_{e_O}$ and in generators $x_e$ with $e\neq e_{O'}$ for
every $2$-cell $O'$, and $x_{e_O}$ occurs in it exactly once, with exponent
$\pm1$. Hence $r_O$ determines a word $w_O$ with $x_{e_O}=w_O$, and $w_O$
contains no generator of the form $x_{e_{O'}}$. Replacing $x_{e_O}$ by $w_O$ in
all relators, then deleting $x_{e_O}$ and $r_O$, is a Tietze transformation. We
check that it is legitimate. Either the parent edge of $O^*$ joins $O^*$ to
$\partial^*$, or it joins $O^*$ to $O_0^*$, in which case the two sides of $e_O$
are $O$ and $O_0$. In the first case $x_{e_O}$
occurs in no relator other than $r_O$, so deleting $x_{e_O}$ and $r_O$ leaves
every remaining relator unchanged. In the second case $x_{e_O}$ occurs in the
relator of another $2$-cell only for $O_0$. There the substitution leaves every
relator other than $r_{O_0}$ unchanged, and the new relator of $O_0$ still
contains $x_{e_{O_0}}$ exactly once, because $w_O$ contains no generator of the
form $x_{e_{O'}}$. The process can therefore
be continued, and it terminates because $\mathbb T$ is finite. At the end every
relator and every generator $x_{e_O}$ has been deleted, and we are left with the
presentation
\[
 \pi_1(Y)\cong\langle\,B\mid\;\rangle,\qquad
 B=\{x_e:\ e\in Y^{(1)}\setminus T,\ e\neq e_O\ \text{for all }O\}.
\]
Hence $B$ is a free basis of $\pi_1(Y)$, computed from the cell-incidence
data alone, and each element is a closed edge-path. Moreover, given any closed edge-path of $Y^{(1)}$ based at the basepoint,
these substitutions compute a word in $B^{\pm1}$ representing its class in
$\pi_1(Y)$.

\subsection{A pro-$p$-closed free subgroup of a closed surface}\label{sec:boundary-extract}

\begin{proposition}\label{prop:free-intermediate}
Let $W$ be a closed connected orientable tiled surface of genus $g(W)\geq2$,
let $Y\subseteq W$ be a compact, connected, boundary reduced tiled subsurface
with nonempty boundary, write $i:Y\hookrightarrow W$ for the inclusion, and fix
a basepoint $w_0$ that is a vertex of $Y^{(1)}$. Then there is a free group $L$
with $i_*(\pi_1(Y,w_0))\leq L\leq\pi_1(W,w_0)$ such that
\begin{enumerate}[label=\textup{(\arabic*)}]
\item $L$ has a computable free basis $B$ whose elements are represented by
      closed edge-paths of $W$;
\item given any closed edge-path of $Y$ based at $w_0$, a word in $B^{\pm1}$
      representing the corresponding element of $i_*(\pi_1(Y,w_0))$ is
      computable;
\item $L$ is pro-$p$-closed in $\pi_1(W,w_0)$, and its pro-$p$ topology is
      induced from $\pi_1(W,w_0)$.
\end{enumerate}
\end{proposition}

\begin{proof}
If $\pi_1(Y,w_0)=1$, take $L=1$. Its free basis is empty, it is pro-$p$-closed
in $\pi_1(W,w_0)$, and its pro-$p$ topology is the one induced from
$\pi_1(W,w_0)$. Henceforth assume $\pi_1(Y,w_0)\neq1$.

Write $N=N(Y)$ for the thick version of $Y$ (see Definition~\ref{def:thick}).
It is a small closed regular neighbourhood of $Y$, it does not depend on the
ambient surface, and we take it in $W$.  Let $\alpha$ be a component of $\partial N$. 
Since $N$ is homotopy equivalent to $Y$ and $\pi_1(Y,w_0)\neq1$, $N$ is not
a disk. Therefore $\alpha$ is nonnull-homotopic in $N$. By Proposition~\ref{prop:pi1-injective}, the inclusion $i$ is $\pi_1$-injective because $Y$ is boundary reduced. Since $Y\hookrightarrow N$ is a
homotopy equivalence, $N\hookrightarrow W$ is $\pi_1$-injective as well. Hence $\alpha$ is an essential simple closed curve in $W$. Moreover, $\alpha$ can be read
off from the combinatorial data of $Y$ by Proposition~3.4 of \cite{mp-core}.

Cutting $W$ along $\alpha$ yields a compact 
 surface $W_\alpha$, together with a
quotient map $q:W_\alpha\to W$ that identifies the boundary circles coming from
the two sides of $\alpha$. Since $Y$ lies
in the interior of $N$, cutting along $\alpha$ does not split $Y$. So there
is a subsurface $Y^0\subseteq W_\alpha$ such that $q$ restricts to a
homeomorphism $q^0:=q|_{Y^0}:Y^0\to Y$. Choose $w_0'\in Y^0$ with 
$q(w_0')=w_0$, and base $\pi_1(W_\alpha)$ at $w_0'$. For the inclusion
$i:Y\hookrightarrow W$ we then have $i\circ q^0=q\circ j$, where
$j:Y^0\hookrightarrow W_\alpha$ is the inclusion. Passing to fundamental groups
gives the commutative diagram
\[
\begin{CD}
 \pi_1(Y^0,w_0') @>{q^0_*}>> \pi_1(Y,w_0) \\
 @V{j_*}VV                  @VV{i_*}V \\
 \pi_1(W_{\alpha},w_0') @>{q_*}>> \pi_1(W,w_0)
\end{CD}
\]
where $q^0_*:\pi_1(Y^0,w_0')\to\pi_1(Y,w_0)$ is an isomorphism, and $i_*$ is
injective because $Y$ is boundary reduced.
Let $W_{\alpha,1}\subseteq W_\alpha$ be the connected component of $W_\alpha$
containing $Y^0$, and let $L=q_*(\pi_1(W_{\alpha,1},w_0'))$. We have
\[
 i_*(\pi_1(Y,w_0))\;=\;q_*\circ j_*\circ(q^0_*)^{-1}(\pi_1(Y,w_0))
 \;\subseteq\; L.
\]

The type
of $\alpha$ is determined by its class in $H_1(W;\Ftwo)$. It is
\emph{separating} when this class is $0$.

\emph{Case 1: $\alpha$ is nonseparating.}
Then $W_\alpha$ is connected, so $W_{\alpha,1}=W_\alpha$ is a compact surface of
genus $g(W)-1$ with two boundary circles, and $\pi_1(W_{\alpha,1})$ is free of rank
$2g(W)-1$. Reattaching $\alpha^+$ to $\alpha^-$ recovers $W$ and is the HNN
extension $\pi_1(W)\cong\pi_1(W_{\alpha,1})\ast_{\langle\alpha\rangle}$; its
base-group inclusion is injective, so $q_*$ is injective on
$\pi_1(W_{\alpha,1},w_0')$. Hence $L$
is a free subgroup of $\pi_1(W,w_0)$ of rank $2g(W)-1$.

\emph{Case 2: $\alpha$ is separating.}
Then $W_\alpha$ has two components, each a connected compact surface with a
single boundary circle; the component $W_{\alpha,1}$ containing $Y^0$ has genus
$h$, say, and the other component is written $W_{\alpha,2}$. As shown above, $\alpha$ is essential in $W$.
Hence neither component of $W_\alpha$ is a disk, and
$1\leq h\leq g(W)-1$. Reattaching
$\alpha^+$ to $\alpha^-$ recovers $W$ and is the free product with amalgamation
$\pi_1(W)\cong\pi_1(W_{\alpha,1})\ast_{\langle\alpha\rangle}\pi_1(W_{\alpha,2})$,
whose base-group inclusions are injective; so $q_*$ is injective on
$\pi_1(W_{\alpha,1},w_0')$. Hence $L$
is a free subgroup of $\pi_1(W,w_0)$ of rank $2h$.

(1) A free basis $B$ of $L$ is computed as follows. Subdivide $W$ so that $Y$ and $\alpha$ are
subcomplexes, and take the cells to be triangles. Cutting the subdivided surface
along $\alpha$ and restricting to the component $W_{\alpha,1}$ produces a finite
cell structure on $W_{\alpha,1}$, together with a recorded cellular gluing map
$q|_{W_{\alpha,1}}$; each $2$-cell boundary path meets each signed edge at most
once. By Section~\ref{sec:freegen}, this cell structure yields a free basis of
$\pi_1(W_{\alpha,1},w_0')$, each basis element given as an explicit edge-path.
Mapping these paths to $W$ by $q|_{W_{\alpha,1}}$, and applying a finite
cellular homotopy in the subdivided polygons, gives based loops in the original
$1$-skeleton of $W$. Reading the original edge labels gives words for a free
basis $B$ of $L$, each represented by a closed edge-path of $W$. All subdivisions,
cutting operations and edge-path replacements are performed combinatorially,
with the original basepoint $w_0$ fixed and the images of the resulting paths
recorded.

(2) By the diagram above, $i_*(\pi_1(Y,w_0))=q_*(j_*(\pi_1(Y^0,w_0')))$, and
$j_*(\pi_1(Y^0,w_0'))\le\pi_1(W_{\alpha,1},w_0')$. Given a closed edge-path of
$Y$ based at $w_0$, read it in the subdivision of (1); as $Y$ and $\alpha$ are
disjoint subcomplexes, $(q^0)^{-1}$ lifts it to a closed edge-path of $Y^0$
based at $w_0'$. Section~\ref{sec:freegen} computes the word in the free basis
of $\pi_1(W_{\alpha,1},w_0')$, and $q_*$ carries it to a word in $B^{\pm1}$
representing the given element of $i_*(\pi_1(Y,w_0))$.

(3) Since $\alpha$ does not pass through $w_0$, its image under any
self-homeomorphism fixing $w_0$ does not either.

If $\alpha$ is nonseparating, let $\gamma$ be a parallel copy
of the standard meridian $a_1$, closed inside a small disc around $w_0$ by an
arc chosen in the sector between the two half-edges of $a_1$. The standard
generators meet only at $w_0$, so $\gamma$ is disjoint from $a_1$, from
$w_0$, and from $a_2,b_2,\dots,a_{g(W)},b_{g(W)}$. Since $\alpha$ and
$\gamma$ are both nonseparating, some self-homeomorphism $f$ of $W$ carries
$\alpha$ to $\gamma$; the self-homeomorphisms preserving $\gamma$ act
transitively on $W\setminus\gamma$, so we may take $f$ to fix $w_0$. Set
$\varphi:=f_*\in\operatorname{Aut}(\pi_1(W,w_0))$. The loop that follows a
path from $w_0$ to $\gamma$, goes once around $\gamma$, and returns is
homotopic to $a_1$ rel $w_0$. Cutting $W$ open along $\gamma$, the component
containing $w_0$ has fundamental group
$$
 \varphi(L)=\langle\,a_1,a_2,b_2,\dots,a_{g(W)},b_{g(W)}\,\rangle.
$$

If $\alpha$ is separating, let $\gamma$ be the standard
separating curve realizing the split $(h,g(W)-h)$, chosen disjoint from
$w_0$ and with $w_0$ on the side of genus $h$; then $a_1,b_1,\dots,a_h,b_h$
may be chosen on that side and disjoint from $\gamma$. Since $\alpha$ and
$\gamma$ have the same split, some self-homeomorphism $f$ of $W$ carries
$\alpha$ to $\gamma$; it carries the genus-$h$ side of $\alpha$, which
contains $Y$, to the genus-$h$ side of $\gamma$, which contains $w_0$, and
the self-homeomorphisms preserving $\gamma$ and that side act transitively
on its interior, so we may take $f$ to fix $w_0$. Set
$\varphi:=f_*\in\operatorname{Aut}(\pi_1(W,w_0))$. Cutting $W$ open along
$\gamma$, the component containing $w_0$ has fundamental group
$$
 \varphi(L)=\langle\,a_1,b_1,\dots,a_h,b_h\,\rangle.
$$

In both cases $\varphi(L)=\langle S\rangle$ for a proper subset $S$ of a
standard generating set of $\pi_1(W,w_0)$. By
Corollary~\ref{cor:automorphism-closed}, $L$ is pro-$p$-closed in
$\pi_1(W,w_0)$, and its pro-$p$ topology is induced from $\pi_1(W,w_0)$.
\end{proof}

\section{Computing the pro-$p$ closure}\label{sec:algorithm}
The machinery of the preceding sections now becomes an algorithm for the
pro-$p$ closure of a finitely generated subgroup of $\Gamma_g$. The algorithm
runs in two phases. Phase~I is a descent along a chain of finite-index
subgroups, driven by mod-$p$ abelianisations; it stops either when the subgroup
has become dense, or when its index passes a computable bound. In the second
case Phase~II applies. The bound guarantees, through the growing process, a
surface with boundary inside the covering, and the problem is then reduced to a
free group. Section~\ref{sec:growing} supplies the growing process and the
bound; Section~\ref{sec:the-algorithm} gives the algorithm and its correctness.

\subsection{The growing process}\label{sec:growing}

The engine of the algorithm is the \emph{growing process} on tiled surfaces of
Magee--Puder~\cite{mp-cover}. Let $\chi_0$ be a fixed integer, and let $f:Y\to Z$ be a morphism of connected tiled
surfaces, where $Y$ is compact and $Z$ has no boundary. The growing process
begins with the subsurface $W_0=f(Y)$. At each stage one annexes every
$4g$-gon of $Z$ that touches the current subsurface along a boundary edge in
one of the configurations of \cite[Definition~2.10]{mp-cover}. The process
stops once the current subsurface is strongly boundary reduced, or boundary
reduced with Euler characteristic less than $\chi_0$. It is deterministic and
terminates after finitely many steps \cite[Lemma~2.11]{mp-cover}, outputting a
compact, connected, boundary reduced tiled subsurface $W_f\subseteq Z$
containing $f(Y)$. We call $W_f$ the \emph{grown surface} attached to $f$. It
fits into the factorisation
\begin{equation}\label{eq:grow-factor}
 Y\xrightarrow{\ f\ }W_f\xrightarrow{\ \iota_f\ }Z,
\end{equation}
where $\iota_f:W_f\hookrightarrow Z$ is the inclusion.

\begin{corollary}\label{cor:grow-pi1}
Suppose $Y$ is boundary reduced and fix a basepoint $w\in Y$. Then the maps
$$
 \pi_1(Y,w)\xrightarrow{\ f_*\ }\pi_1(W_f,f(w))\xrightarrow{\ (\iota_f)_*\ }\pi_1(Z,f(w))
$$
induced by the factorisation~\eqref{eq:grow-factor} are injective; identifying
each group with its image, we have
$$
 \pi_1(Y,w)\;\leq\;\pi_1(W_f,f(w))\;\leq\;\pi_1(Z,f(w)).
$$
\end{corollary}
\begin{proof}
Since $Y$ and $W_f$ are boundary reduced, Proposition~\ref{prop:pi1-injective}
shows that the morphism $f:Y\to W_f$ and the inclusion $\iota_f:W_f\hookrightarrow Z$
are both $\pi_1$-injective. The induced maps on fundamental groups are therefore
injective, and we identify each group with its image in the next.
\end{proof}

The growth of the process is bounded by the starting surface $Y$ alone.
By Lemma~2.12 of \cite{mp-cover} and its proof, the process annexes at
most $B(Y)$ $4g$-gons in total. Here $B(Y)\geq 0$ is a computable bound depending
only on the combinatorial data of $Y$, $g$ and $\chi_0$. In particular it is independent
of $f$ and of $Z$. The initial subsurface $W_0=f(Y)$ already has
at most $\#\mathrm{polygons}(Y)$ $2$-cells, so $W_f$ has no more than
\[
 \widehat B(Y):=\#\mathrm{polygons}(Y)+B(Y)
\]
$4g$-gons. Hence $\widehat B(Y)$ is computable from the combinatorial data of
$Y$ and the fixed $\chi_0$ alone.

\begin{corollary}\label{cor:no-fill}
Let $K\leq\Gamma_g$ be a finite-index subgroup and $Y$ a compact connected
tiled surface. Let $f:Y\to\Core(K)$ be a morphism of tiled surfaces. If
$[\Gamma_g:K]>\widehat B(Y)$ for the fixed $\chi_0$, then $W_f$ is a proper,
connected, boundary reduced subsurface of $\Core(K)$ with nonempty boundary.
\end{corollary}
\begin{proof}
Since $K$ has finite index, $\Upsilon_K$ is a compact finite cover of
$\Sigma_g$; every complement component in Definition~\ref{def:core} then
contains finitely many $2$-cells, so all are adjoined and $\Core(K)=\Upsilon_K$.
It is closed and has exactly $[\Gamma_g:K]$ $4g$-gons, one per sheet (since
$\Sigma_g$ has a single $2$-cell, its $n$-sheeted cover carries exactly $n$
$2$-cells). By the bound $\widehat B(Y)$ above, $W_f$ has fewer $4g$-gons
than $\Core(K)$, so $W_f\ne\Core(K)$. A proper subsurface of a closed surface
has nonempty boundary, and $W_f$ is boundary reduced by construction.
\end{proof}

\subsection{The algorithm}\label{sec:the-algorithm}
Throughout this subsection, for $J\le\Gamma_g$ let $\Upsilon_J\to\Sigma_g$ be
the covering surface corresponding to $J$ and let $\widetilde o_J\in\Upsilon_J$
be the lift of $o$; if $H\le J$, the covering morphism $\Upsilon_H\to\Upsilon_J$
sends $\widetilde o_H$ to $\widetilde o_J$. The core surface $\Core(J)$ need not
contain $\widetilde o_J$, so we first pass to a conjugate subgroup for which it
does; this preserves the isomorphism type of the core surface as a tiled surface.

\begin{lemma}\label{lem:standard-form}
Let $1\neq H\leq\Gamma_g$ be given by words $w_1,\dots,w_m$ in
$A_0^{\pm1}$. One can compute $s\in\Gamma_g$ such that, after replacing $H$ by the conjugate
$sHs^{-1}$ and each $w_j$ by $s\,w_j\,s^{-1}$, the subgroup (again denoted
$H$) satisfies $\widetilde o_H\in\Core(H)$. This replacement preserves the isomorphism
type of $\Core(H)$ as a tiled surface and leaves
$\widehat B(\Core(H))$ unchanged.
\end{lemma}
\begin{proof}
Dehn's algorithm \cite{lyndon-schupp} discards the $w_j$ representing $1$;
since $H\ne1$, at least one word remains. Keep all remaining words and
relabel one of them $w_1$. Let $C$ be the closed path in $(\Sigma_g)^{(1)}$
reading $w_1$. The surface $\Sigma_g$ has no boundary cycles, hence is a
boundary reduced tiled surface, so by Corollary~4.10 of \cite{mp-core} the
free homotopy class of $C$ has a shortest representative cycle in
$(\Sigma_g)^{(1)}$, reached by the finitely many length-decreasing steps of
that proof. Reading it gives a word $u$ of least length in the conjugacy class
of $w_1$. Enumerate the words
in $A_0^{\pm1}$ by length and, for each of them, say $s$, test the equation
$us=sw_1$, which Dehn's algorithm decides; the enumeration meets a conjugator
$s$ with $u=s\,w_1\,s^{-1}$ and stops there.

Then $u=s\,w_1\,s^{-1}$ lies in $sHs^{-1}$, so the path in
$\Upsilon_{sHs^{-1}}$ reading $u$ from $\widetilde o_{sHs^{-1}}$ is closed.
It has length $|u|$ and is a shortest representative cycle of its free
homotopy class in $\Upsilon_{sHs^{-1}}$, since a shorter cycle would project
to a shorter word in $\Gamma_g$ conjugate to $u$. By
Definition~\ref{def:core} the core surface contains every shortest
representative cycle, so $\widetilde o_{sHs^{-1}}\in\Core(sHs^{-1})$.

Conjugate subgroups have isomorphic covering surfaces, via an isomorphism
induced by a deck transformation of the universal cover; commuting with the
projection to $\Sigma_g$, it preserves labels and directions, hence carries
$\Core(H)$ onto $\Core(sHs^{-1})$ and preserves the combinatorial data
defining $\widehat B$. Hence the replacement preserves the isomorphism type
of $\Core(H)$ as a tiled surface, and it leaves $\widehat B(\Core(H))$
unchanged.
\end{proof}

\begin{lemma}\label{lem:covering-morphism}
Let $1\ne H\le K\le\Gamma_g$ with $K$ of finite index and $\widetilde o_H\in\Core(H)$.
There is a unique morphism of tiled surfaces $f:\Core(H)\to\Core(K)$ with
$f(\widetilde o_H)=\widetilde o_K$. The covering projections to $\Sigma_g$
identify $\pi_1(\Core(H),\widetilde o_H)$ and $\pi_1(\Core(K),\widetilde o_K)$
with the subgroups $H$ and $K$ of $\Gamma_g=\pi_1(\Sigma_g,o)$, and under
these identifications $f_*$ is the inclusion $H\le K$. Moreover, $f$ can be
computed from $\Core(H)$ and the right action of $\Gamma_g$ on the set
$K\backslash\Gamma_g$ of right cosets.
\end{lemma}
\begin{proof}
Since $K$ has finite index, $\Core(K)=\Upsilon_K$, as in the proof of
Corollary~\ref{cor:no-fill}. The vertices of
$\Upsilon_J$ are the lifts of $o$. Fix one of them, $\widetilde o_J$. Every
vertex is the endpoint of a path from $\widetilde o_J$ reading some word
$g\in\Gamma_g$, and that endpoint is $\widetilde o_J$ exactly when $g\in J$.
Two words $g,g'$ lead to the same vertex exactly when $gg'^{-1}\in J$, so
reading words from $\widetilde o_J$ identifies the vertices with the right
cosets $J\backslash\Gamma_g=\{Jg:g\in\Gamma_g\}$, with $\widetilde o_J$
corresponding to $J$; reading a further word $r$ from $Jg$ ends at $Jgr$.

Under this identification the covering morphism $\Upsilon_H\to\Upsilon_K$
induced by $H\le K$ carries $Hg$ to $Kg$. This is well defined because
$H\le K$, and it sends $\widetilde o_H$, the coset $H$, to the coset
$K=\widetilde o_K$.
Its restriction to the subcomplex $\Core(H)\subseteq\Upsilon_H$ is a morphism
$f:\Core(H)\to\Core(K)$ with $f(\widetilde o_H)=\widetilde o_K$. A morphism of
tiled surfaces preserves labels and directions, and at each vertex of
$\Upsilon_K$ there is at most one incoming and one outgoing edge of a given
label. So the image of a vertex determines the images of all edges, and
property \textup{(P3)} then determines the $2$-cells. As $\Core(H)$ is
connected, $f$ is the unique morphism with $f(\widetilde o_H)=\widetilde o_K$.

The covering projections $p_H:\Core(H)\to\Sigma_g$ and
$p_K:\Core(K)\to\Sigma_g$ satisfy $p_K\circ f=p_H$, because $f$ is the
restriction of a covering morphism. Hence the diagram
\[
\begin{CD}
\pi_1(\Core(H),\widetilde o_H) @>{f_*}>> \pi_1(\Core(K),\widetilde o_K)\\
@V{(p_H)_*}VV @VV{(p_K)_*}V\\
H @>{\subseteq}>> K
\end{CD}
\]
commutes, and the two vertical maps are isomorphisms onto $H$ and onto $K$. For
the left one this is Proposition~\ref{prop:mp-basic}, applied to the embedding
$\Core(H)\hookrightarrow\Upsilon_H$ (here $H\ne1$); the right one is the
covering projection of $\Upsilon_K=\Core(K)$. Hence $f_*$ is the inclusion
$H\le K$.

To compute $f$, fix a spanning tree $T$ of $\Core(H)$ rooted at
$\widetilde o_H$, and for each vertex $v$ let $r_v\in A_0^{\pm1}$ be the word
read along the tree path from $\widetilde o_H$ to $v$. Put $f(v):=Kr_v$, the
vertex reached from $\widetilde o_K$ by reading $r_v$. Each edge of $\Core(H)$
leaves a vertex $v$ with a label and a direction, and $\Upsilon_K$ has at most
one edge leaving $f(v)$ with this label and direction, so the image of every
edge is forced.
Now let $r_v'$ be another tree word for $v$. Both words reach $v$, so
$Hr_v=Hr_v'$ and hence $r_v=hr_v'$ with $h\in H\le K$. Reading $r_v$ and
$r_v'$ from $\widetilde o_K$ gives $Kr_v=Khr_v'=Kr_v'$, because $Kh=K$, so
$f(v)$ does not depend on $T$.
\end{proof}

The algorithm below works with a descending chain of finite-index subgroups $K_i$ of
$\Gamma_g$, beginning at $K_0=\Gamma_g$. At the start of round $i$ it holds a
generating set $A_i$ of $K_i$, each element written, after recursive
substitution, as a word in $A_0^{\pm1}$, a presentation of $K_i$ in these
generators;
for $i=0$ these are the standard data of $\Gamma_g$. Write $\sigma_i:K_i\to
(\Fp)^{2g(K_i)}$ for the mod-$p$ abelianisation, $g(K_i)$ being the genus of $\Upsilon_{K_i}$. Since $K_{i+1}=\sigma_i^{-1}(\sigma_i(H))$ contains $H$, we have $H\le K_i$
throughout. Each $K_i$ is clopen in $\Gamma_g$. This holds for $K_0=\Gamma_g$. If $K_i$
is clopen, then its pro-$p$ topology is induced from $\Gamma_g$
(Proposition~\ref{prop:clopen-topology}), while $K_{i+1}$ is normal in $K_i$
with quotient  the finite $p$-group $Q_i$; hence $K_{i+1}$ is a
basic open set of the pro-$p$ topology of $K_i$, so it is open, and therefore
closed, in $\Gamma_g$. The algorithm runs as follows.

\begin{algorithm}[H]
\small
\setlength{\abovedisplayskip}{4pt plus 1pt}
\setlength{\belowdisplayskip}{4pt plus 1pt}
\caption{Computing the pro-$p$ closure of a finitely generated subgroup}
\label{alg:main}
\begin{algorithmic}
\Require A prime $p$, an integer $g\geq2$, the alphabet
       $A_0=\{a_1,b_1,\dots,a_g,b_g\}$, and words $w_1,\dots,w_m$ in
       $A_0^{\pm1}$, $m\ge1$, generating $H$.
\Ensure A finite set of words in $A_0^{\pm1}$ generating $\Cl_p(H)$.
\Statex \textbf{Preprocessing.}
\Statex Free-reduce the input words; this does not change $H$, so we may
        assume the $w_j$ freely reduced from now on. If every $w_j$ represents
        $1$ (word problem), output the empty set and
        stop. Otherwise apply Lemma~\ref{lem:standard-form}, which renames $H$
        and the $w_j$ and returns $s$; keep $s$ for the output. Build
        $\Core(H)$ as in Proposition~\ref{prop:mp-basic}, and compute
        $\widehat B(\Core(H))$ for a fixed integer $\chi_0$, say
        $\chi_0=2$. Set $K_0:=\Gamma_g$
        and $i:=0$.
\Statex \textbf{Phase~I: density descent.}
\State The following three steps are repeated, starting from the current
       $i$.
\State (1) If $[\Gamma_g:K_i]>\widehat B(\Core(H))$, proceed to Phase~II.
\State (2) Compute $\sigma_i(H)$, the $\Fp$-span of the images of the words
       $w_1,\dots,w_m$ under $\sigma_i$.
\State (3) If $\sigma_i(H)=(\Fp)^{2g(K_i)}$, output $s^{-1}A_i s$ and stop.
       Otherwise $K_{i+1}:=\sigma_i^{-1}(\sigma_i(H))$ is the kernel of the
       composite epimorphism
$$
 \ K_i\xrightarrow{\ \sigma_i\ }(\Fp)^{2g(K_i)}
 \twoheadrightarrow Q_i:=(\Fp)^{2g(K_i)}\big/\sigma_i(H),
$$
where $Q_i$ is a finite elementary abelian $p$-group of order $p^{d_i}$
with $d_i=2g(K_i)-\dim_{\Fp}\sigma_i(H)\ge1$, so the index is updated by
$[\Gamma_g:K_{i+1}]=p^{d_i}[\Gamma_g:K_i]$. By Reidemeister--Schreier,
applied with any right transversal of $K_{i+1}$ in $K_i$ (one is obtained
from the finite quotient $Q_i$), the presentation of $K_i$ in $A_i$ yields a
generating set $A_{i+1}$ of $K_{i+1}$, its elements written as words in
$A_i^{\pm1}$, together with a presentation of $K_{i+1}$ in these generators and a rule that carries each word in $A_i^{\pm1}$
representing an element of $K_{i+1}$ to the word in $A_{i+1}^{\pm1}$
representing the same element. The map $\sigma_{i+1}$ is computed from the
presentation of $K_{i+1}$ in $A_{i+1}$ by the linear algebra over
$\mathbb Z$. Rewrite
$w_1,\dots,w_m$ by this rule, replace $i$ by $i+1$ and repeat from~(1).

\Statex \textbf{Phase~II: reduction to a free group.}
\State (1) By coset enumeration on the known index, compute the action of
       $\Gamma_g$ on the right cosets $K_i\backslash\Gamma_g$; this yields the
       whole cover $\Core(K_i)=\Upsilon_{K_i}$, together with the morphism
       $f_i:\Core(H)\to\Core(K_i)$ of
       Lemma~\ref{lem:covering-morphism}, so that
       $f_i(\widetilde o_H)=\widetilde o_{K_i}$.
\State (2) By step~(1) of Phase~I, $[\Gamma_g:K_i]>\widehat B(\Core(H))$,
       so Corollary~\ref{cor:no-fill} applies. Run the growing process with
       this same $\chi_0$ to obtain the boundary reduced subsurface
       $W_{f_i}\subseteq\Core(K_i)$ with nonempty boundary. Since $\widetilde o_H\in\Core(H)$, the
       point $\widetilde o_{K_i}=f_i(\widetilde o_H)$ is a vertex of
       $f_i(\Core(H))\subseteq W_{f_i}$.  By
       Corollary~\ref{cor:grow-pi1}, $H=\pi_1(\Core(H),\widetilde o_H)$ embeds
       in $\pi_1(W_{f_i},\widetilde o_{K_i})$.
\State (3) Apply Proposition~\ref{prop:free-intermediate} to $W=\Core(K_i)$,
       $Y=W_{f_i}$ and $w_0=\widetilde o_{K_i}$.  It yields a free subgroup
       $L\le K_i=\pi_1(\Core(K_i),\widetilde o_{K_i})$ with
$$
 H=\pi_1(\Core(H),\widetilde o_H)\;\le\;
 \pi_1(W_{f_i},\widetilde o_{K_i})\;\le\;L.
$$
       The subgroup $L$ is pro-$p$-closed in $K_i$ with the pro-$p$ topology
       induced from $K_i$, and it has a computable free basis $B$ whose
       elements are closed edge-paths of $\Core(K_i)$, hence words in
       $A_0^{\pm1}$.
\State (4) Apply the free-group algorithm \cite{msw01,rzk-add,weil98} to $H$
       in the free group $L$ with free basis $B$, the generators of $H$ being
       written in $B^{\pm1}$ by
       Proposition~\ref{prop:free-intermediate}\textup{(2)}. It returns words
       $v_1,\dots,v_r$ in $B^{\pm1}$ generating $\Cl_p(H,L)$. Replace each
       letter of $B$ by its word in $A_0^{\pm1}$ and output
       $s^{-1}v_1s,\dots,s^{-1}v_rs$.

\end{algorithmic}
\end{algorithm}

\begin{theorem}\label{thm:termination-correctness}
Algorithm~\ref{alg:main} terminates after finitely many steps and outputs a
generating set of $\Cl_p(H)$.
\end{theorem}
\begin{proof}

If $H=1$, the preprocessing step outputs the empty set, which generates
$\Cl_p(1)=1$. Hence assume $H\ne1$. Lemma~\ref{lem:standard-form} produces
$s\in\Gamma_g$ such that, after renaming, $\widetilde o_H\in\Core(H)$; the
original subgroup is $s^{-1}Hs$. Conjugation by $s$ is an automorphism of
$\Gamma_g$, hence a homeomorphism of its pro-$p$ topology (the argument is
that of Lemma~\ref{lem:automorphism-invariance}), so
\[
 \Cl_p(s^{-1}Hs,\Gamma_g)=s^{-1}\,\Cl_p(H,\Gamma_g)\,s .
\]
It suffices to prove the claim for the renamed subgroup; we return to the
conjugation at the end.

{\bf Termination.}
The construction in \cite[Theorem~5.10]{mp-core} begins with a preparation
phase; Lemma~\ref{lem:standard-form} supplies its first step, conjugating the
generating set so that $w_1$ becomes its shortest conjugate. The construction
then performs the remaining preparatory steps and the finite folding process,
which builds $\Core(H)$. The recorded image of the wedge point is
$\widetilde o_H$, which lies in $\Core(H)$ by
Lemma~\ref{lem:standard-form}, as Lemma~\ref{lem:covering-morphism} requires.
The bound $\widehat B(\Core(H))$ is computed from the combinatorial data of
$\Core(H)$.

In Phase~I, suppose the algorithm reaches round $i$ without exiting. The
homomorphism $\sigma_i$ is determined by its values on the generators $A_i$;
these values are computed from the presentation of $K_i$ by linear algebra
over $\mathbb Z$. Hence each $\sigma_i(w_j)$, and therefore the span
$\sigma_i(H)$, is computable, and the density test is a rank computation. If
the test fails,
$d_i\ge1$, and $K_{i+1}$, the kernel of an epimorphism onto $Q_i$ of order
$p^{d_i}$, satisfies
\[
 [\Gamma_g:K_{i+1}]=p^{d_i}\,[\Gamma_g:K_i]\ \ge\ p\,[\Gamma_g:K_i].
\]
A descent step occurs only while $[\Gamma_g:K_i]\le\widehat B(\Core(H))$;
since each step multiplies the index by at least $p$, a $k$-th step requires
$p^{k-1}\le\widehat B(\Core(H))$, so at most
$\lfloor\log_p\widehat B(\Core(H))\rfloor+1$ steps occur (and none when
$\widehat B(\Core(H))=0$). Phase~II is therefore entered only with
$[\Gamma_g:K_i]>\widehat B(\Core(H))$.

Assume Phase~II is entered.

(1) Coset enumeration on the known index gives
$\Core(K_i)=\Upsilon_{K_i}$, and Lemma~\ref{lem:covering-morphism} computes
$f_i$.

(2) By Phase~I, $[\Gamma_g:K_i]>\widehat B(\Core(H))$; Corollary~\ref{cor:no-fill}
then gives a proper boundary reduced $W_{f_i}$ with nonempty boundary, and the
growing process terminates \cite[Lemma~2.11]{mp-cover}.

(3) Proposition~\ref{prop:free-intermediate} and Section~\ref{sec:freegen}
perform finitely many combinatorial operations, so this step is finite.

(4) A spanning tree of $\Core(H)^{(1)}$ rooted at $\widetilde o_H$, together
with its non-tree edges, gives closed paths $\beta_e$ generating
$H=\pi_1(\Core(H),\widetilde o_H)$. Their images $f_i(\beta_e)$ are closed
edge-paths of $W_{f_i}$ based at $\widetilde o_{K_i}$, and
Proposition~\ref{prop:free-intermediate}\textup{(2)} writes them in $B^{\pm1}$;
this is the finite generating set required by step~(4) of Phase~II. The
free-group algorithm of \cite{msw01,rzk-add,weil98} then runs in finitely many
steps on this input.

{\bf Correctness.}
If Phase~I halts at the density test, then
$\sigma_i(H)=(\Fp)^{2g(K_i)}$. As $K_i$ has finite index in $\Gamma_g$, it is
a surface group with $H_1(K_i)\cong\mathbb Z^{2g(K_i)}$, so $\sigma_i$ is its
mod-$p$ abelianisation and Proposition~\ref{prop:density-equiv} applies
verbatim, giving $H$ pro-$p$ dense in $K_i$. Since $K_i$ is clopen and
$\mathbf G_p$ is extension closed, its pro-$p$ topology is induced from
$\Gamma_g$ (Proposition~\ref{prop:clopen-topology}); hence
$\Cl_p(H,\Gamma_g)\supseteq K_i$, and the reverse inclusion holds because
$K_i$ is pro-$p$-closed and contains $H$. Thus $\Cl_p(H,\Gamma_g)=K_i$ is
generated by $A_i$.

If Phase~II runs to completion, identify $H$ with
$\pi_1(\Core(H),\widetilde o_H)$ by
Proposition~\ref{prop:mp-basic}\textup{(1)}. Then
\[
 H=\pi_1(\Core(H),\widetilde o_H)\;\le\;\pi_1(W_{f_i},\widetilde o_{K_i})\;\le\;L\;\le\;K_i .
\]
The first inclusion is Corollary~\ref{cor:grow-pi1}, applied to the
factorisation~\eqref{eq:grow-factor} for $f_i$: $\Core(H)$ is connected and
boundary reduced, $\widetilde o_H$ is a vertex, and
$f_i(\widetilde o_H)=\widetilde o_{K_i}$. Proposition~\ref{prop:free-intermediate},
applied as in step~(3) of Phase~II to $W=\Core(K_i)$ (of genus
$g(K_i)\ge2$), $Y=W_{f_i}$ and $w_0=\widetilde o_{K_i}$ (a vertex of
$W_{f_i}^{(1)}$, as noted in step~(2)), gives the second and third inclusions
and shows that $L$ is free, pro-$p$-closed in $K_i$, with its pro-$p$ topology
induced from $K_i$. As $K_i$ is clopen in $\Gamma_g$, its topology is induced
from $\Gamma_g$ (Proposition~\ref{prop:clopen-topology}) and it is closed
there. Proposition~\ref{prop:free-reduction} applied to
$H\le L\le K_i\le\Gamma_g$ now gives
\[
 \Cl_p(H,L)=\Cl_p(H,K_i)=\Cl_p(H,\Gamma_g).
\]
The free-group algorithm, run on the words of item~(4), returns
$v_1,\dots,v_r\in B^{\pm1}$ generating $\Cl_p(H,L)$; replacing each letter of
$B$ by its word in $A_0^{\pm1}$
(Proposition~\ref{prop:free-intermediate}\textup{(1)}) yields words in
$A_0^{\pm1}$ generating $\Cl_p(H,\Gamma_g)$.

In both exits the algorithm outputs a finite set of words in $A_0^{\pm1}$
generating $\Cl_p(H,\Gamma_g)$; conjugating by $s^{-1}$ gives generators of
$s^{-1}\Cl_p(H,\Gamma_g)s=\Cl_p(s^{-1}Hs,\Gamma_g)$, the closure of the
original input.
\end{proof}

\section*{Tool Use}

All mathematical content is due to the authors. GPT-6.1 Sol and DeepSeek-V4.1-Flash were used to comment on the draft and to draft and polish English sentences, including the wording of some definitions, statements and proofs. The authors are responsible for all the contents of this paper.

\end{document}